\documentclass[11pt,reqno]{amsart}
\usepackage{amsthm, amsmath}
\usepackage{siunitx}
\usepackage{color}
\usepackage{mathrsfs}
\newcommand {\Q}  {{\mathbb Q}}
\newcommand {\Z}  {{\mathbb Z}}

\newcommand {\C}  {{\mathbb C}}

\newcommand {\PP}  {{\mathbb P}}
\newcommand {\N}  {{\mathbb N}}

\newcommand {\g}   {\mathfrak g}
\newcommand {\HH}  {{\mathcal H}}
\newcommand{\ord}{\textup{ord}}
\newcommand{\m}{\textup{m}}

\newcommand{\uu} {u}

\newtheorem {Th}{Theorem}

\newtheorem {Pro}{Proposition}
\newtheorem {Le}{Lemma}

\begin{document}
\title{Composite polynomials in linear recurrence sequences}
\author{Clemens Fuchs, Christina Karolus}
\address[Clemens Fuchs]{University of Salzburg, Hellbrunnerstr. 34/I, 5020 Salzburg, Austria}
\email{clemens.fuchs@sbg.ac.at}
\address[Christina Karolus]{University of Salzburg, Hellbrunnerstr. 34/I, 5020 Salzburg, Austria}
\email{christina.karolus@sbg.ac.at}
\keywords{decomposable polynomials, linear recurrence sequences, Brownawell-Masser inequality}
\subjclass[2010]{}

\begin{abstract}
Let  $(G_n(x))_{n=0}^\infty$ be a $d$-th order linear recurrence sequence having polynomial characteristic roots, one of which has degree strictly greater than the others. Moreover, let $m\geq 2$ be a given integer. We ask for $n\in\N$ such that the equation $G_n(x)=g\circ h$ is satisfied for a polynomial $g\in\C[x]$ with $\deg g=m$ and some polynomial $h\in\C[x]$ with $\deg h>1$. We prove that for all but finitely many $n$ these decompositions can be described in ``finite terms'' coming from a generic decomposition parameterized by an algebraic variety. All data in this description will be shown to be effectively computable.
%We prove that if $G_n$ is decomposable in the requested form infinitely often, then the set of such indices $n\in\N$ is the union of a finite set and finitely many arithmetic progressions. For the finite set, we give an explicit upper bound. It follows that, whenever $G_n(x)=g\circ h(x)$ holds, $h$ is an element of one out of finitely many polynomial recurrence sequences. Moreover, we give a complete description of possible recurrence sequences $(H_n)_{n=0}^\infty$ such that $h\in (H_n)$.}
\end{abstract}
\maketitle

\section{Introduction and results}
\label{intro}

\let\thefootnote\relax\footnote{Supported by Austrian Science Fund (FWF) Grant No. P24574.}

Let $\C[x]$ be the polynomial ring in the indeterminate $x$ (we remark right away that $\C$ might be replaced by an algebraically closed field of characteristic $0$; however, all polynomials below are assumed to have coefficients in $\C$). Composition of polynomials is a well-defined operation on $\C[x]$. It is associative and has with $f(x)=x$ an identity element, but it is neither commutative nor distributive. There are many reasons to be interested in polynomial composition, e.g. if for a given $f\in\C[x]$ there is a $g\in\C[x]$ with $f\circ g$ irreducible over $\C$, then $f$ is an irreducible polynomial in $\C[x]$ as well. Another example is that the decompositions of $f$ exhibit arithmetical properties associated with $f$, which is used to solve equations of separated variable type (cf. \cite{BT,Bea}). The invertible elements in $\C[x]$ with respect to decomposition are the linear polynomials. We call $f(x)=g\circ h$ a non-trivial decomposition if neither $g$ nor $h$ is linear. Let $m\geq 2$ be an integer; we call $f(x)=g\circ h$ an $m$-decomposition if $\deg g=m$ and we say that $f$ is $m$-decomposable if an $m$-decomposition of $f$ exists. We call $f$ indecomposable if $f$ admits only trivial decompositions. A pair $(g,h)\in\C[x]$ is called equivalent to $(g',h')$ if there are $a,b\in\C,a\neq 0$ such that $g(x)=g'(ax+b),h(x)=(h'(x)-b)/a$. It is easy to see that every polynomial $f\in\C[x]$ can be decomposed as $f(x)=f_1\circ f_2\circ\cdots \circ f_k$ with $f_i$ indecomposable. Moreover, this decomposition is unique in the following sense: if $f(x)=f_1'\circ\cdots\circ f_l'$ with $f_j'$ indecomposable is another decomposition, then $k=l$ and $f_1,\ldots,f_k$ are obtained by replacing neighboring pairs of $f_1',\ldots,f_l'$ by an equivalent one a finite number of times.  (This is known as Ritt's first theorem; cf. \cite{S}.) There is a nice algebraic description of decompositions of a polynomial $f\in\C[x]$ since they are (up to equivalence) in one-to-one correspondence to intermediate fields between $\C(x)$ and $\C(f(x))$ (cf. again \cite{S}).

We start with a few general remarks. First, it might happen that for a given $f\in\C[x]$ and a given $g\in\C[x]$ we have different $h_1,h_2\in\C[x]$ such that $f(x)=g\circ h_1=g\circ h_2$. However, in this case $g(X)=g(Y)$ has a solution $X=h_1(x),Y=h_2(x)$. This situation was completely solved in \cite{AZ}. It follows that $g=g'\circ x^k$ for some $k>1$ and $h_1$ and $h_2$ just differ by a constant (to be more precise, by a $k$-th root of unity). Second, by linear equivalence we can control the leading coefficient and the constant term of $h$. E.g. we may assume that $h\in\C[x]$ is monic and satisfies $h(0)=0$ (every other value in $\C$ is fine as well). Moreover, assume that $f$ is given and that we have $f(x)=g\circ h$ with $f,g,h\in\C[x]$. Let $a\in\C\backslash\{0\}$ be the leading coefficient of $f$. Then we may also assume that $g$ is monic by writing $f(x)=a(g\circ h)$. Third, assume that $f,h\in\C[x]$ are given. Then there are at most finitely many $g\in\C[x]$ with $f(x)=g\circ h$. This can be seen as follows: We may assume that $f,g,h$ are all monic. Write $g(x)=(x-b_1)\cdots(x-b_m)$, where $b_1,\ldots,b_m\in\C$ are not necessarily distinct. Assume that $f(x)=(x-a_1)\cdots(x-a_n)$ with $a_1,\ldots,a_n\in\C$. Then $g(h(x))=(h(x)-b_1)\cdots(h(x)-b_m)=(x-a_1)\cdots(x-a_n)=f(x)$. It follows that there is a partition of the multi-set $\{a_1,\ldots,a_n\}$ with equally large blocks (of size $\deg h$) that describe $g$ uniquely. If we assume that $h(0)=0$, then the $b_i$ are just the product of all elements in the $i$-th block. The unique $g$ can be found, without calculating the roots of $f$, by comparing coefficients in $f(x)=g\circ h$ (cf. \cite{R}).

In this paper we are interested in non-trivial decompositions (with two factors, an ``inner'' and an ``outer'' factor) of polynomials with coefficients in $\C$. This problem is hard in general since the decompositions of polynomials can be anything \emph{a priori} (since conversely, every pair $(g,h)$ gives a polynomial $g\circ h$). Therefore, it is natural to restrict to a subset of $\C[x]$ which is described by a finite amount of data and then to ask whether or not all decompositions in this subset can be described \emph{in finite terms} depending on the data describing the subset. We give a few (important and non-trivial) examples to illustrate this approach.

Let $n\geq 2$ be a given integer. We consider the set of all polynomials $f\in\C[x]$ of degree $n$. Then there is an integer $J$ and for every $1\leq j\leq J$ an algebraic variety $\mathcal{V}_j\subset \mathbb{A}^{n+t_j}$ for some $2\leq t_j\leq n$ defined over $\Q$ for which equations can be written down effectively and there are polynomials $f_j,h_j,g_j$ with coefficients in the coordinate ring of the variety and depending on integers $k_1,\ldots,k_{t_j}$ and $(l_1,\ldots,l_n)\in\{0,1,\ldots,n-1\}^n$ such that the following holds: a) $g_j\circ h_j=f_j$ is a polynomial of degree $n$ with coefficients in the coordinate ring; b) for every point $P\in\mathcal{V}_j(\C)$ and integers $k_1,\ldots,k_{t_j},l_1,\ldots,l_n$  one gets a decomposition $f_j(P,x)=g_j(P,h_j(P,x))$; c) conversely, for every polynomial $f\in\C[x]$ of degree $n$ and every non-trivial decomposition $f(x)=g\circ h$ with $g(x)$ not of the shape $(ax+b)^m,m\in\N,a,b\in\C$ there are $P\in\mathcal{V}_j(\C),k_1,\ldots,k_{t_j},l_1,\ldots,l_n$ such that $f(x)=f_j(P,x),g(x)=g_j(P,x),h(x)=h_j(P,x)$. This result (formulated in different forms) can be found in \cite{FP2,B,R}.

Let $\ell$ be a given integer. We consider the set of \emph{lacunary} polynomials (with respect to $\ell$), that is the set of all polynomials $f\in\C[x]$ with $\ell$ non-constant terms. Then there are integers $p,J$ depending on $\ell$ and for every $1\leq j\leq J$ an algebraic variety $\mathcal{V}_j$ defined over $\Q$ and a lattice $\Lambda_j$ for which equations can be written down explicitly and (Laurent-)polynomials $f_j,h_j\in\Q[\mathcal{V}_j][z_1^{\pm 1},\ldots,z_p^{\pm 1}],g_j\in\Q[\mathcal{V}_j][z]$ with coefficients in the coordinate ring of the variety such that the following holds: a) $g_j\circ h_j=f_j$ is a (Laurent-)polynomial with $\ell$ non-constant terms with coefficients in the coordinate ring; b) for every point $P\in\mathcal{V}_j(\C)$ and $(u_1,\ldots,u_p)\in\Lambda_j$ one gets a decomposition $f_j(P,x^{u_1},\ldots,x^{u_p})=g_j(P,h_j(P,x^{u_1},\ldots,x^{u_p}))$; c) conversely, for every polynomial $f\in\C[x]$ with $\ell$ non-constant terms and every non-trivial decomposition $f(x)=g\circ h$ with $h(x)$ not of the shape $ax^m+b,m\in\N,a,b\in\C$ there is a $j$, a point $P\in\mathcal{V}_j(\C)$ and $(u_1,\ldots,u_p)\in\Lambda_j$ such that $f(x)=f_j(P,x^{u_1},\ldots,x^{u_p}),g(x)=g_j(P,x),h(x)=h_j(P,x^{u_1},\ldots,x^{u_p})$. This result can be found in \cite{Z1}; cf. also \cite{Z3,Z4}. (A similar result holds for lacunary rational functions $f\in\C(x)$ by a combination of \cite{FZ} and \cite{FMZ}.)

In the present paper we are interested in another subset of $\C[x]$ namely the subset $\{G_n(x); n\in\N\}$ that consists of elements of a linear recurrence sequence $(G_n(x))_{n=0}^\infty$ of polynomials in $\C[x]$. The sequence is fixed by the recurrence relation and by the initial values. Equivalently, every element of the sequence can be written by a Binet-type formula $G_n(x)=a_1\alpha_1^n+\cdots+a_t\alpha_t^n$, where $\alpha_1,\ldots,\alpha_t$ are the distinct roots of the characteristic polynomial associated to the recurring relation and the $a_i$ are polynomials in $n$ with coefficients in the splitting field $\C(x,\alpha_1,\ldots,\alpha_t)$ of degree less than the corresponding multiplicity of $\alpha_i$ as a root of the characteristic equation. In this way all elements are given by a finite amount of data. Our goal is to describe all decomposable $G_n$'s in this set and all their decompositions in finite terms, depending only on the given data.
To fix terms we shall consider the $d$-th order linear recurrence sequence $(G_n(x))_{n=0}^\infty$, given by the relation
\begin{align}\label{recrel}
G_{n+d}(x)=A_{d-1}(x)G_{n+d-1}(x)+\cdots+A_0(x)G_n(x),
\end{align}
with $A_0,\ldots,A_{d-1}\in\C\left[x\right]$ and initial terms $G_0,\ldots,G_{d-1}\in\C[x]$. Denote by $\alpha_1,\ldots,\alpha_t$ the distinct characteristic roots of the sequence, that is the characteristic polynomial $\mathcal{G}\in\C(x)[T]$ splits as
\[
\mathcal{G}(T)=T^d-A_{d-1}T^{d-1}-\cdots - A_0=(T-\alpha_1)^{k_1}(T-\alpha_2)^{k_2}\cdots(T-\alpha_t)^{k_t},
\]
where $k_1,\ldots,k_t\in\N$. We assume that all roots are simple, i.e. $t=d$, and that they are polynomials, i.e. $\alpha_i\in\C[x]$ for $i=1,\ldots,d$. Then $G_n(x)$ admits a representation of the form
\begin{align}
G_n(x)=a_1\alpha_1^n+a_2\alpha_2^n+\cdots+a_d\alpha_d^n.\label{binet}
\end{align}
By assumption we consider the special situation that $a_1,\ldots,a_d\in\C$ and $\alpha_1,\ldots,\alpha_t\in\C[x]$. Finally, we assume that $\deg(\alpha_1)>\deg(\alpha_i)$ for $i>1$.
% Can we rephrase this condition in an equivalent condition for $G_i, A_i,\ i=0,\ldots,d-1$?

We mention that for binary recurrences the authors together with Kreso proved in \cite{FKK} that if $G_n(x)=g\circ h$, then either $\deg g$ is bounded independently of $n$ and only in terms of the initial data, unless $h$ is special (meaning that $(g,h)$ is equivalent to $(g',x^m)$ or $(g'',T_m(x))$ where $(T_n(x))_{n=0}^\infty$ denotes the sequence of Chebyshev polynomials and $g',g''\in\C[x]$) or a technical condition is not verified (see the paper for details). This describes the ``outer'' decomposition factor in such a decomposition. In view of this result, which we expect (without the technical condition) to hold in general, we restrict ourselves to $m$-decompositions for an integer $m\geq 2$ which we view as fixed from now on.

We further mention that for a given sequence $(G_n(x))_{n=0}^\infty$ the decompositions of the form $G_n(x)=G_m\circ h$ for a fixed polynomial $h\in\C[x], \deg h\geq 2$ were considered in \cite{FPT1,F,FPT2}. It was Zannier who proved in general that this equation has only finitely many solutions $(n,m),n\neq m$, unless we are in the cyclic or Chebyshev case as above (cf. \cite{Z2}). This result was made effective in \cite{FP1}. A further result in this direction can be found in \cite{FPT3}.

There are a few trivial situations that we have to take into account below. If $G_n(x)=f(\beta^n)$ with $f,\beta\in\C[x]$, then every decomposition $f(x)=g\circ h$ with $\deg g=m$ leads to a sought decomposition $G_n(x)=g(h(\beta^n))$ for every $n\in\N$. Observe that this situation might also lead to slightly different decompositions. Assume e.g. that $G_n(x)=a_1\alpha_1^n+a_2,a_1,a_2\in\C,\alpha_1\in\C[x]$; if $n$ is a multiple of $m$, i.e. $n=m\ell$, then $G_n(x)=g\circ h$ with $g(x)=a_1x^m+a_2,h(x)=\alpha_1^\ell$. More generally, when $G_n(x)=g(H_n(x))$ with $g\in\C[x],\deg g=m$ and $(H_n(x))_{n=0}^\infty$ is another linear recurrence sequence in $\C[x]$, then obviously we again have a sought decomposition for every $n\in\N$. Unfortunately it seems that these cases are not exhaustive. There might be many ``sporadic'' solutions that arise by polynomial-exponential equations that are complicated to control in general.
%Polynomial decomposition: definition, out vs. inner decomposition factor, equivalence of decompositions
%Decompositions of a single polynomial, of all polynomials with given degree, of lacunary polynomials, of polynomials with given number of zeroes
%Setting
%Resultate erwähnen: Zannier (lakunäre Polynome); ev. andere Lakunaritäten (Fuchs-Pethö); Rekursionen über Funktionenkörper (Zannier, Fuchs-Pethö-Tichy); Avanzi-Zannier
%Offensichtliche Familien: $G_n(x)=a_1\alpha_1^n+a_2,a_1,a_2\in\C,\alpha_1\in\C[x]$; if $n$ is a multiple of $m$, i.e. $n=ml$, then $G_n(x)=g\circ h$ with $g(x)=a_1x^m+a_2,h(x)=\alpha_1^l$
%Fuchs-Karolus-Kreso
%Equivalence implies that we can normalize the leading coefficient and the constant term: $h$ should be taken monic (i.e. $g$ has leading coefficient $=a_1$) and with $h(0)=0$
%Motivating examples:

We start with the following theorem, which clarifies the structure of the ``inner'' decomposition factor that may appear in an $m$-decomposition of elements in the sequence $(G_n(x))_{n=0}^\infty$.

\begin{Th}\label{th1}
Let $(G_n(x))_{n=0}^\infty$ be a non-degenerate simple linear recurrence sequence of order $d\geq 2$ with power sum representation $G_n(x)=a_1\alpha_1^n+\cdots+a_d\alpha_d^n$ with $a_1,\ldots,a_d\in\C$, $\alpha_1,\ldots,\alpha_d\in\C[x]$ satisfying $\deg\alpha_1>\max\{\deg\alpha_2,\ldots,\deg\alpha_d\}$. Moreover, let $m\geq 2$ be an integer. Write $m_0$ for the least integer such that $\alpha_1^{m_0/m}\in \C[x]$. Then there is an effectively computable positive constant $C$ such that the following holds: Assume that for some $n\in\N$ with $n>C$ we have $G_n(x)=g\circ h$ with $\deg g=m,\deg h>1$. Then there are $c_1,\ldots,c_l\in\C$ such that \[h(x)=c_1\gamma_1^\ell+\cdots +c_l\gamma_l^\ell,\] where $m_0\ell=n$ and $l\in\N$ is bounded explicitly in terms of $m,d$ and $\deg(\alpha_1)+\cdots+\deg(\alpha_d)$ and $\gamma_1,\ldots,\gamma_l\in\C(x)$ can be given explicitly in terms of $\alpha_1,\ldots,\alpha_d$, both independently of $n$.
\end{Th}

This result should be compared with Proposition 2 in \cite{Z1}.

We illustrate the result with an example. Let $(G_n(x))_{n=0}^\infty$ be given by $G_n(x)=x^{3n}+3(2x^2)^n+3(4x)^n+2^{3n}$ for all $n\geq 0$ and let $m=3$. We have $\alpha_1=x^3,\alpha_2=2x^2,\alpha_3=4x,\alpha_4=8$ and $m_0=1$. The proof of the theorem shows that we must have $n<10440$ or $h(x)=c_1x^n+c_2$ with $c_1,c_2\in\C$. Let $g(x)=x^3$. Then $g(h(x))=(c_1x^n+c_2)^3=c_1^3x^{3n}+3c_1^2c_2x^{2n}+3c_1c_2^2x^n+c_2^3$. Comparing $g(h(x))$ with $G_n(x)$ shows that $c_1^3=1,c_1^2c_2=2^n,c_1c_2^2=4^n,c_2^3=8^n$. This defines a subvariety $\mathcal{V}$ of $\mathbb{A}^2\times \mathbb{G}_\m$. Up to (possibly) finitely many exceptions for small $n$ we have $G_n(x)=g(h(x))=(c_1x^n+c_2)^3$, where $(c_1,c_2,n)\in\mathcal{V}(\C)$.

Observe that $m_0$ in the theorem can also be described as follows: Write $\alpha_1(x)=v(x-v_1)^{k_1}\cdots(x-v_t)^{k_t}$ and define $\psi$ by $\psi^m=\alpha_1$. Put $d=$ gcd$(k_1,\ldots,k_t,m)$. Then $m_0=m/d$. Obviously, $\psi^{m/d}\in\C[x]$. Conversely, observe that $m_0$ is a divisor of $m$ since by definition of $m_0$ the polynomial $T^{m_0}-\psi^{m_0}$ is the minimal polynomial of $\psi$ over $\C(x)$ (cf. Proposition \ref{kummer}) and thus divides $T^m-\psi^m$ over $\C(x,\psi)$. Since $\psi^{m_0}\in\C(x)$, it follows $m_0k_i/m\in\N$ and thus $m/m_0$ divides $k_i$ (thus gcd$(k_1,\ldots,k_t,m)=d$) for $i=1,\ldots,t$. The smallest such integer is obtained in the case of equality giving $m=m_0d$ as claimed.

The structure of all $m$-decompositions for a given $m\geq 2$ can now be described as follows.

\begin{Th}\label{th2}
Let $(G_n(x))_{n=0}^\infty$ be a non-degenerate simple linear recurrence sequence of order $d\geq 2$ with power sum representation $G_n(x)=a_1\alpha_1^n+\cdots+a_d\alpha_d^n$ with $a_1,\ldots,a_d\in\C$, $\alpha_1,\ldots,\alpha_d\in\C[x]$ satisfying $\deg\alpha_1>\max\{\deg\alpha_2,\ldots,\deg\alpha_d\}$. Moreover, let $m\geq 2$ be an integer. Write $m_0$ for the least integer such that $\alpha_1^{m_0/m}\in \C[x]$. Then there is an explicitly computable positive constant $C$, and a subvariety $\mathcal{V}$ of $\mathbb{A}^{l+m+1}\times\mathbb{G}_\m^{t}$ with $t,l$ bounded explicitly in terms of $m,d$ and $\deg(\alpha_1)+\cdots+\deg(\alpha_d)$ for which a system of polynomial-exponential equations in the polynomial variables $c_1,\ldots,c_l,g_0,\ldots,g_m$ and the exponential variable $\ell$ (with coefficients in $\Q$) can be written down explicitly such that the following holds:
\begin{itemize}
\item[a)] Defining $G(x)=g_0x^m+g_1x^{m-1}+\cdots +g_m\in\C[\mathcal{V}][x]$ and $H_\ell=c_1\gamma_1^\ell+c_2\gamma_2^\ell+\cdots +c_l\gamma_l^\ell\in\C[\mathcal{V}][x]$, where $\gamma_1,\ldots,\gamma_l\in\C(x)$ can be given explicitly in terms of $\alpha_1,\ldots,\alpha_d$, then $G_{m_0\ell}=G\circ H_\ell$ holds as an equation in $x$ with coefficients in the coordinate ring of $\mathcal{V}$. In particular, for any point $P=(c_1,\ldots,c_l,g_0,\ldots,g_m,\ell)\in\mathcal{V}(\C)$ we get a decomposition $G_n(x)=g\circ h$,  $g(x)=G(P,x)\in\C[x]$ and $h(x)=H_l(P,x)\in\C[x]$ (with $n=m_0\ell$).
\item[b)] Conversely, let $G_n(x)=g\circ h$ be a decomposition of $G_n(x)$ for some $n\in\N$ with $g,h\in\C[x],\deg g=m, \deg h>1$. Then either $n\leq C$ or there exists a point $P=(g_0,\ldots,g_m,c_1,\ldots,c_l,\ell)\in\mathcal{V}(\C)$ with $g(x)=G(P,x)$ and $h(x)=H_\ell(P,x)$ and $n=m_0\ell$.
\end{itemize}
\end{Th}

\noindent Remarks and special cases:\\[0.2cm]
{ a)} Binary case: Let $(G_n(x))_{n=0}^\infty$ be a non-degenerate binary simple linear recurrence which does not satisfy a recurrence relation of order less than $2$; thus, we have $G_n(x)=a_1\alpha_1^n+a_2\alpha_2^n$ with $a_1,a_2\in\C$. We assume that $\alpha_1,\alpha_2\in\C[x]$ and $\deg\alpha_1>\deg\alpha_2$. Moreover, we assume that one of the conditions of \cite[Theorem 2]{FKK} is satisfied. Then there is an effectively computable constant $C$ and there are finitely many subvarieties $\mathcal{V}_i$ of $\mathbb{A}^{m_i+1+l_i}\times\mathbb{G}_\m^{t_i}$ and equations $G_{m_{0,i}\ell}(x)=G^{(i)}\circ H^{(i)}_\ell$, where $\deg G^{(i)}=m_i\geq 2$, in the coordinate ring of $\mathcal{V}_i$ such that the following holds: If $G_n(x)=g\circ h$ for some $n\in\N$ and $g,h\in\C[x],\deg g,\deg h>1$ with $h(x)$ indecomposable and not of the shape $ax^m+b,m\in\N,a,b\in\C$, then either $n\leq C$ or there is an $i$ and a $P=(g_{i0},\ldots,g_{im_i},c_{i1},\ldots,c_{il_i},\ell)\in\mathcal{V}_i(\C)$ such that $n=m_{0,i}\ell, g(x)=G^{(i)}(P,x),h(x)=H^{(i)}_\ell(P,x)$.\\[0.2cm]
{ b)} When all $\alpha_i$ are monic, then the varieties can be chosen without the $\mathbb{G}_\m$-part.\\[0.2cm]
{ c)} Assume that $\alpha_1\in\C[x]$ satisfies $1\le k:=\deg\alpha_1\le m$. %and $\alpha_2\cdots\alpha_d$ do not have common roots.
Then $G_n(x)=g\circ h$ with $n>C$ implies that $h$ is of the form $c_1+c_2\alpha_1^{\deg G_n/(mk)}$ with $c_1,c_2\in\C$.\\[0.2cm]
%up to leading coefficient not all powers of a given polynomial
{ d)} Assume that $\alpha_1=\beta^{m_1},\alpha_2=\beta^{m_2},\ldots,\alpha_d=\beta^{m_d}$ with $m_1>m_2\geq \cdots\geq m_d\geq 0$. Then $G_n(x)=f(\beta^n)$, where $f(x)=f_1x^{m_1}+f_2x^{m_2}+\cdots +f_dx^{m_d}$. In this case it follows that either $n\leq C$ or $h(x)=c_1\beta^{k_1\ell}+\cdots+c_l\beta^{k_l\ell}=h'(\beta^\ell)$ for some $h'\in\C[x]$. Thus if $G_n(x)=g\circ h=(g\circ h')(\beta^\ell)=f(\beta^{m_0\ell})$. Therefore, the problem reduces to find all $m$-decompositions of the polynomial $f\circ x^{m_0}$.\\[0.2cm]
{ e)} The proof shows that if $(G_n(x))_{n=0}^\infty$ is defined over a number field, i.e. the coefficients of the Binet-type equation (\ref{binet}) as well as the characteristic roots are polynomials with coefficients in some number field $K$, then all decomposition factors $g,h$ are defined over $K$ as well. In this case we are interested in decompositions over $K$, which can be described by the above conclusion of the statements.\\[0.2cm]
%\item Whenever we are given with a recurrence $(G_n(x))_{n=0}^\infty$ defined over a number field $K$ and such that $h(x)=c_1\gamma_1^\ell+c_2\gamma_2^\ell+c_3\gamma_3^\ell$ with $c_1,c_2,c_3\in K$, then it is possible (by using [Corvaja-Zannier] to prove that for any given $g\in K[x]$ of degree $m\geq 2$ up to finitely many exception there are at most finitely many infinite families of solutions, which can be calculated explicitly, that can be parametrized by linear recurrence sequences (i.e. we get a decomposition $G_{a\ell+b}(x)=g(H_\ell)$ for all $\ell\in\N$). A concrete example is given by $G_n(x)=x^{6n}+(x^4(x-1))^n$, where we get $h(x)=c_1x^{2n}+c_2+c_3(x-1)^n$.
{ f)} We also remark that the above results include a description in finite terms of all $m$-th powers in a linear recurring sequence of polynomials satisfying the conditions of the theorem (i.e. the sequence is non-degenerate and simple and the characteristic roots are polynomials where one has degree larger than all others). This follows by fixing $g(x)=x^m$ and then going through the proof of the above theorem.\\[0.2cm]
{ g)} Finally, we mention that if we know that the $c_i$ can be parametrized by power sums as well (in particular if they are constant) and that we have a decomposition for any $\ell$ (or for all members along an arithmetic progression), then these families are easy to calculate. This follows since we may identify varying powers by indeterminates (see e.g. \cite[Lemma 2.1]{CZ}) and then use the algorithm in \cite{B} for polynomials in several variables (actually, we view such a polynomial as a polynomial in one of the variables; the other variables can be embedded into $\C$ so that we may view the polynomial again as an element in $\C[x]$) to determine the decompositions.\\

The proof of the theorems follows essentially the ideas of \cite{Z1}. Assume that $G_n(x)=g\circ h$. This equality is viewed as an equation for the unknown $h=h(x)$; it is a root of $g(T)-G_n(x)=0$ over the (rational) function field $\C(G_n(x))$. Thus we can expand $h$ as a Puiseux series in terms of quantities $G_n(x)^{s/m}, s=1,0,-1,\ldots$, where $m=\deg g$. Then one uses the multinomial series to expand $G_n(x)^{s/m}$ for any $s$; in order to justify this multiple expansion, the ``dominant root condition'' on the degrees of the characteristic roots is needed. Afterwards we use, as in \cite{Z1}, a function field variant of the Schmidt subspace theorem (Proposition \ref{zannier}) proved in \cite{Z1}, to find that either $n$ is bounded or $h$ can be expressed as given in Theorem \ref{th1}. Using this information, one views the $c_1,\ldots,c_l$ as well as the coefficients $g_0,\ldots,g_m$ of $g$, while the degree of $g$ is fixed, as indeterminants and then compares $g\circ h$ with $G_n(x)$ for the given $n\in\N$. Using unit equations over function fields, this either implies that $n$ is bounded or we have two linear recurrences that are related (see \cite{SS} for this notation). In the latter case, up to a permutation, the characteristic roots have to match up and then, since they are monic polynomials at that point, the coefficients coincide. This gives polynomial-exponential equations that can be written down explicitly and which define a variety. From this the statement follows.

The rest of the paper is organized as follows. In the next section we collect some auxiliary results that will be needed for the proof of the theorems. In Section 3 we give the proofs of Theorem \ref{th1} and \ref{th2}. In Section 4 we give some more details justifying the remarks and special cases.

\section{Auxiliary results}
\label{sec:2}

In this section, we recall some basic information and collect some statements, which we will make use of in our proofs later on.

Generally, an algebraic function field $F/K$ is a finite algebraic extension of $K(x)$, where $x$ is some element transcendental over $K$. If $F$ is itself of the shape $F=K(x)$, then $F$ is said to be rational. The rational function field has genus $\g_F=0$. Throughout this paper, we will work over the complex numbers $K=\C$, even though our proofs hold over any other algebraically closed field as well. Then
\[
\C(x)=\left\{\frac{f(x)}{g(x)};\ f(x),g(x)\in\C[x],\ g(x)\neq 0\right\},
\]
i.e. $\C(x)$ is the field of fractions of $\C[x]$. On $\C(x)$ we define valuations as follows. For each $a\in\C$, let $\nu_a(f)$ be the unique integer such that $f(x)=(x-a)^{\nu_a(f)}p(x)/q(x)$, where $p,q\in\C[x]$ are such that $p(a)q(a)\neq0$. Moreover, with the symbol $\infty$ we associate the valuation $\nu_{\infty}(f)=\deg q-\deg p$, where $f(x)=p(x)/q(x)$. If $\nu_a(f)>0$ for an $a\in\C$, $a$ is called a zero of $f$, and it is called a pole of $f$, if $\nu_a(f)<0$. These functions are all (normalized, up to equivalence) valuations on $\C(x)$ and for a finite extension $L$ of $\C(x)$ each one of them can be extended to at most $[F:\C(x)]$  valuations on $F$, which again gives all discrete valuations on $F$. Both, in $\C(x)$ and in $F$, for any $f\in F/\C$ the so-called sum formula holds, that is
\[
\sum\limits_{\nu}\nu(f)=0,
\]
where the sum is taken over all valuations on the respective function field. There is a one-to-one relation between valuations and places, namely for any valuation
$\nu_a$ on $\C(x)$, $a\in\C\cup\{\infty\}$,  there is a place
$P_a=\{f\in\C(x);\ \nu_a(f)>0\}$
(it is the unique maximal ideal of the valuation ring
$\mathcal{O}_{P_a}=\{f\in\C(x);\ \nu_a(f)\ge0\}$).
Therefore, valuations are sometimes introduced in terms of places (and often instead of $P_a$ we simply write $a$). We write $\mathbb{P}_F$ for the set of places of the field $F$. Now let $F'$ be an algebraic extension of the function field $F$. Then a place $P'\in\mathbb{P}_{F'}$ is said to lie over $P\in\mathbb{P}_F$, if $P\subset P'$. We write $P'|P$ in this case. Then there is an integer $e=e(P'|P)$,  $1\le e\le[F':F]$, called the ramification index of $P'$ over $P$, such that $\nu_{P'}(x)=e\cdot\nu_P(x)$ for all $x\in F$. We say that $P'|P$ is ramified if $e(P'|P)>1$ and unramified otherwise. If $e(P'|P)=[F':F]$, there is exactly one place $P'\in F'$ lying above $P\in F$ and $P'$ is said to be totally ramified. The places $P'\in\mathbb{P}_{F'}$ lying above $P\in\mathbb{P}_F$ correspond to the extensions of the respective valuation $\nu_P$ in $F$.
Denote by $F_P$ the residue class field $\mathcal{O}_{P}/P$. Then the relative degree of $P'$ over $P$  is defined to be $f(P'|P)=[F'_{P'}:F_P]$.
We shall need the following two statements, which can be found in \cite{St}.

\begin{Pro}\label{eisenstein}
Let $F/\C$ be a function field in one variable and $\varphi\in F\left[T\right]$,
\begin{align*}
\varphi(T)=a_nT^n+a_{n-1}T^{n-1}+\cdots+a_1T+a_0.
\end{align*}
If there is a place $P\in\PP_F$ such that $\nu_P(a_n)=0$, $\nu_P(a_i)\ge0$ for $i=1,\ldots,n-1$, $\nu_P(a_0)<0$ and $\gcd(n,\nu_P(a_0))=1$, then $\varphi(T)$ is irreducible in $F\left[T\right]$. Furthermore, if $F'=F(y)$, where $y$ is a root of $\varphi(T)$, then $P$ has a unique extension $P'\in \PP_{F'}$, $e(P'|P)=n$ and $f(P'|P)=1$.
\end{Pro}

\begin{Pro}\label{kummer}
Let $F/\C$ be a function field in one variable. Suppose that $u\in F$ satisfies $u\neq w^d$ for all $w\in F$ and $d|n$, $d>1$. Let $F'=F(z)$ with $z^n=u$. Then $F'$ is said to be a Kummer extension of $F$ and we have:
\begin{itemize}
\item[(a)] The polynomial $\varphi(T)=T^n-u$ is the minimal polynomial of $z$ over $F$ (in particular, it is irreducible over $F$). The extension $F'/F$ is Galois of degree $n$; its Galois group is cyclic, and all automorphisms of $F'/F$ are given by $\sigma(z)=\zeta z$, where $\zeta\in\C$ is an $n$-th root of unity.
\item[(b)] Let $P\in\mathbb{P}_F$ and $P'\in\mathbb{P}_{F'}$ be an extension of $P$. Let $r_P:=\gcd(n, \nu_P(u))$. Then $e(P'|P)=n/r_P$.
\item[(c)]Denote by $g$ (resp. $g'$) the genus of $F/\C$ (resp. $F'/\C$). Then
\[
g'=1+n(g-1)+\frac{1}{2}\sum\limits_{P\in\mathbb{P}_F}(n-r_P)\deg P.
\]
\end{itemize}
\end{Pro}

Our strategy also involves the use of height functions in function fields. Define the  projective height $\mathcal{H}$ of $u_1,\ldots,u_n\in F/\C$, where $n\ge 2$ and not all $u_i$ are zero, via
\begin{align*}
\mathcal{H}(u_1,\ldots,u_n)=-\sum\limits_{\nu}\min(\nu(u_1),\ldots,\nu(u_n)).
\end{align*}
Also, for a single element $f\in F^*$, set
\begin{align*}
\mathcal{H}(f)=-\sum\limits_{\nu}\min(0,\nu(f)).
\end{align*}
In both cases the sum is taken over all discrete valuations $\nu$ on $F$. Note that $\nu(f)\neq0$ only for a finite number of valuations $\nu$ and that $\HH(f)=\sum_{\nu}\max (0,\nu(f))$ if $f\in F^*$, by the sum formula. For $f=0$, we define $\HH(f)=\infty$. We state some basic properties of the projective height, cf. \cite{FKK}.

\begin{Le} \label{propert}
Denote as above by $\mathcal{H}$ the projective height on $F/\C$. Then for $f,g\in F^*$ the following properties hold:
\begin{center}
\begin{enumerate}
\item $\HH(f)\ge 0$ and $\HH(f)=\HH(1/f)$,
\item $\HH(f)-\HH(g)\le\mathcal{H}(f+g)\le \HH(f)+\HH(g)$,\label{plus}
\item $\HH(f)-\HH(g)\le\HH(fg)\le \HH(f)+\HH(g)$,\label{mult}
\item $\HH(f^n)=|n|\cdot\HH(f)$,
\item $\HH(f)=0\Leftrightarrow f\in \C^*$,%should we note in the statement that we assume $k=\bar{k}$?
\item $\HH(A(f))=\deg A\cdot \HH(f)$ for any $A\in \C[T]\backslash\left\{0\right\}$.\label{last}
\end{enumerate}
\end{center}
\end{Le}

The following proposition is an important ingredient for the proof of our first theorem. It can be seen as a function field analogue of the Schmidt subspace theorem, modelled by Zannier, cf. \cite{Z1}.

\begin{Pro}[Zannier]\label{zannier}
Let $F/\C$ be a function field in one variable, of genus $\g$, let $\varphi_1,\ldots,\varphi_n\in F$ be linearly independent over $\C$ and let $r\in\{0,1,\ldots,n\}$. Let $S$ be a finite set of places of $F$ containing all the poles of $\varphi_1,\ldots,\varphi_n$ and all the zeros of $\varphi_1,\ldots,\varphi_r$. Put $\sigma=\sum_{i=1}^n\varphi_i$. Then
\begin{align*}
\sum\limits_{\nu\in S}\left(\nu(\sigma)-\min_{i=1,\ldots,n}\nu(\varphi_i)\right)\le\binom{n}{2}(\vert S\vert+2\g-2)+\sum_{i=r+1}^n\deg(\varphi_i).
\end{align*}
\end{Pro}

Recall that for a finite set $S$ of places (or valuations, respectively) of $F$, an element $f\in F$ is called an $S$-unit, if it has zeros and poles only at places in $S$, i.e. the set of $S$-units is given by
\[
\mathcal{O}_S^*=\{f\in L; \nu(f)=0\mbox{ for all }\nu\notin S\}.
\]
We will also use the following result due to Brownawell and Masser \cite{BM} taken from \cite{FZ}, giving an upper bound for the height of $S$-units, which arise as a solution of certain $S$-unit-equations.

\begin{Pro}[Brownawell-Masser] \label{bm}
Let  $F/\C$ be a function field in one variable of genus $\g$. Moreover, for a finite set $S$ of discrete valuations, let  $u_1,\ldots, u_n$ be $S$-units, not all constant,  and
$$1+u_1+u_2+\cdots +u_n=0,$$
where no proper subsum of the left side vanishes. Then it holds
\begin{align*}
\max\limits_{i=1,\ldots,n}\HH(u_i)\le\frac{1}{2}(n-1)(n-2)(|S|+2\g-2).
\end{align*}
\end{Pro}

In our proof we will use an expansion of the polynomial $h$ in $G_n=g\circ h$ as a Puiseux series. Therefore, we give a quick review on formal power series, Laurent series and Puiseux series (cf. \cite{W} and \cite{Setal}).
Formally, a (complex) polynomial is a sequence $(a_0,a_1,\ldots)$, where $a_i\in\C$ and where there exists an $n\in\N$ such that $a_j=0$ for all $j\ge n$. Such an element is associated with the finite sum $a_0+a_1x+\ldots+a_nx^n$, where $x$ is an indeterminate identified with the element $(0,1,0,0,\ldots)$. Together with the usual addition and multiplication this gives the $\C$-algebra of (complex) polynomials $\C[x]$. If we also allow sequences with infinite support, we obtain the algebra of formal power series, denoted by $\C[\![x]\!]$, that is the set
\[
\C[\![x]\!]=\left\{(a_0,a_1,a_2,\ldots);\ a_i\in\C\right\}=\left\{a_0+a_1x+a_2x^2+\ldots;\ a_i\in\C\right\},
\]
where addition and multiplication are defined just in the same way as for polynomials. Note that the notation as an infinite sum is meant only formally, i.e. questions of convergence are disregarded. $\C[\![x]\!]$ is an integral domain and the units in $\C[\![x]\!]$ are precisely the elements with non-zero constant term.
The quotient ring of $\C[\![x]\!]$ is the ring of formal Laurent series, denoted by $\C(\!(x)\!)$. It is the localization of $\C[\![x]\!]$ with respect to the ideal $(x)$ and its elements are given by the set
\begin{align*}
\C(\!(x)\!)&=\left\{(a_m,a_{m+1},a_{m+2},\ldots);\ m\in\Z,\ a_i\in\C\right\}\\
&=\left\{a_mx^m+a_{m+1}x^{m+1}+a_{m+2}x^{m+2}+\ldots;\ m\in\Z,\ a_i\in\C\right\},
\end{align*}
so a Laurent series is the sum of a formal power series plus possibly a finite number of terms with negative exponent. The ring $\C[x]$ carries the topology inherited from the product topology of $\C^\N$ given by pointwise convergence. $\C[\![x]\!]$ is the completion of the polynomial ring $\C[x]$ with respect to this topology. Defining $\ord(f)$ to be the smallest $i$ such that $a_i\neq0$ if $f\neq 0$, and $\ord(0)=\infty$ gives a discrete valuation on $\C(\!(x)\!)$. The valuation ring is given by $\C[\![x]\!]$ and the residue field is $\C$. Moreover, $\C(\!(x)\!)$ is a field (in general, $K(\!(x)\!)$ is a field if $K$ is a field), which is the completion with respect to this valuation topology of the field $\C(x)$ of rational functions.
As a generalization, the field of formal Puiseux series is obtained by allowing also fractional exponents, i.e. Laurent series in $\C(\!(x^{1/n})\!)$ for some $n\in\N$. More precisely, the field of formal Puiseux series is given by
\[
\C(\!(x^{1/\infty})\!)=\bigcup\limits_{n=0}^\infty\C(\!(x^{1/n})\!).\]

The valuation ord naturally extends to this field and it is complete with respect to the induced topology. The classical Newton-Puiseux theorem shows that $\C(\!(x^{1/\infty})\!)$ is an algebraic closure of the field of formal Laurent series $\C(\!(x)\!)$.

For the expansion of $h$ as a Puiseux series, we rely on the following classical theorem, cf. \cite{E}.

\begin{Pro}[Puiseux's Theorem] \label{puiseux}
Let $F/\C$ be a function field in one variable of degree $\left[F:\C(x)\right]=n$. Then there are $1\leq r\le n$ natural numbers $e_i$ satisfying
\begin{align*}
e_1+e_2+\cdots+e_r=n
\end{align*}
which have the following meaning: The irreducible equation $f(x,y)=0$ satisfied by an arbitrary function $y$ in $F$ has for solutions the $r$ series
\begin{align}
y_i=\sum\limits_{k=\nu_i}^\infty a_{ik}x^{-k/e_i},\quad a_{i\nu_i}\neq0,\quad i=1,2,\ldots,r.
\end{align}
With a primitive $e_i$-th root of unity $\zeta$ form
\begin{align*}
y_{ij}=\sum\limits_k a_{ik}\zeta^{jk}x^{-k/e_i},\quad j=0,\ldots,e_i-1;
\end{align*}
then $f(x,y)$ is identical with
\begin{align}\label{factor}
f(x,y)=\prod\limits_{i,j}(y-y_{ij}).
\end{align}
The coefficients $a_{ik}$ are elements of a finite field extension $F'$ of $F$ and their images under isomorphisms of $F'$ give permutations of the $y_{ij}$ in \eqref{factor}. %The power series have respective radii of convergence $\neq0$.
\end{Pro}

We remark that in the above theorem $r$ is the number of places $P_i|P_\infty$ for the unique infinite place $P_\infty\in\mathbb{P}_{\C(x)}$, where $P_i\in \mathbb{P}_F$. Moreover, $\nu_i=\nu_{P_i}(y)$ is the valuation of $y$ at $P_i$ and $e_i=e(P_i|P_\infty)$ is the ramification index of $P_i$ over $P_\infty$.

Finally, we state the following little lemma that will be useful in the proof.

\begin{Le}\label{lem1}
Let $f\in\C[x]$. Then $f(1/y)\in\C(y)$ with a pole only at $y=0$. Moreover, the order of vanishing at $y=0$ of $f(1/y)$ is equal to $-\deg f$.
\end{Le}

\emph{Proof}. Write $f(x)=f_0(x-a_1)^{k_1}\cdots (x-a_t)^{k_t}$ with $k_1+\cdots+k_t=\deg f$. Put $y=1/x$. Then $f(1/y)=f(x)=y^{-(k_1+\cdots+k_t)}f_0(1-a_1y)^{k_2}\cdots(1-a_ty)^{k_t}$, which shows that the multiplicity of $y=0$ as a pole of $f(1/y)$ in $\C(y)$ is equal to $\deg f$. This is the claim. \hfill$\square$\\

\section{Proof of Theorem \ref{th1} and \ref{th2}}
\label{sec:3}

\subsection{Proof of Theorem \ref{th1}}

Let $g\in\C\left[x\right]$ be a polynomial of degree $\deg g=m\ge2$ and $(G_n(x))_{n=0}^\infty$ as given in the theorem.

Let $b_1,\ldots,b_d\in\C$ be the leading coefficients of $\alpha_1,\ldots,\alpha_d$ respectively. We write $\alpha_i^n=b_i^n\beta_i^n$ for $i=1,\ldots,d$ and therefore have \[G_n(x)=a_1b_1^n\beta_1^n+\cdots+a_db_d^n\beta_d^n,\] where $a_1,\ldots,a_d,b_1,\ldots,b_d\in\C$ and $\beta_1,\ldots,\beta_d\in\C[x]$ have leading coefficients equal to $1$ and satisfy $\deg\beta_1>\max\{\deg\beta_2,\ldots,\deg\beta_d\}$.

Let $K$  be the rational function field $K=\C(x)$ and let $z$ be a root of $\varphi(T)=g(T)-x=g_0T^m+\cdots +g_m\in K[T]$. At the infinite place $P=P_\infty$ we have $\nu_\infty(g_0)=\nu_\infty(g_i)=0$ for $i=1,\ldots,m-1$ and $\nu_\infty(g_m)=\nu_\infty(g(0)-x)=-1$. Also, $\gcd(\deg\varphi,\nu_\infty(g_m))=\gcd(m,-1)=1$. Therefore, by Proposition \ref{eisenstein}, $\varphi(T)$ is irreducible over $\C(x)$ and  $P_\infty$ has a unique extension $P'_{\infty'}\in\PP_{L}$, where $L=\C(x)(z)=\C(x,z)$. By Puiseux's Theorem (Theorem \ref{puiseux}), it follows that there is an expansion of $z$ of the form
\begin{align}
z&=\sum\limits_{k=\nu_{\infty'}}^\infty \uu_k(\sqrt[m]{1/x})^k\notag\\
&=\uu_{-1}(\sqrt[m]{1/x})^{-1}+\uu_{0}(\sqrt[m]{1/x})^{0}+\uu_{1}(\sqrt[m]{1/x})^{1}+\uu_{2}(\sqrt[m]{1/x})^{2}+\cdots\notag\\
&=\uu_{-1}x^{1/m}+ \uu_{0}+\uu_{1}x^{-1/m}+\uu_{2}x^{-2/m}+\cdots,
\end{align}
where the $\uu_j\in\C$ depend only on $g$. Note that we have used that $e_i(\infty')=m$ in Theorem \ref{puiseux}, by Proposition \ref{eisenstein} and that $\nu_{\infty'}(z)=-1$. This can be seen as follows. Since $z$ is a root of $\varphi(T)$, we have that $g(z)=x$, hence $\nu_{\infty'}(g_0y^m+\ldots+g_{m-1}y+g(0))=\nu_{\infty'}(x)=m\cdot\nu_\infty(x)=-m$. By integrality we have $\nu_{\infty'}(y)<0$, so using the strict triangle inequality we get
\[m\cdot\nu_{\infty'}(y)=\min\limits_{i=1,\ldots,m}(i\cdot\nu_{\infty'}(y))=\min\limits_{i=1,\ldots,m}(\nu_{\infty'}(g_{m-i}y^i))=-m.\]
This expansion is understood as an equality in the algebraic closure of the ring of formal Laurent series $\C(\!(x)\!)$, which itself is the usual metric completion of $K$ (which has the property that an infinite sum converges if and only if each fixed power of $x$ appears in only finitely many terms). Hence, if $G_n(x)=g\circ h$, substituting $G_n(x)$ for $x$ and $h(x)$ for $z$, this yields an expansion for $h(x)$ of the form
\begin{align}
h(x)=\uu_{-1}G_n(x)^{1/m}+\uu_{0}+\uu_{1}G_n(x)^{-1/m}+\uu_{2}G_n(x)^{-2/m}+\cdots,
\end{align}
for a suitable choice of the $m$-th root of $G_n(x)$, where the coefficients $\uu_i\in\C$ depend only on $g$. As in \cite{Z1}, we expand the different roots $G_n(x)^{s/m}$ for $s\in\{1,0,-1,\ldots\}$ using the multinomial theorem and equation \eqref{binet} to get
\begin{equation}
\begin{split}
G_n&(x)^{s/m}\\
&=(a_1\alpha_1^n+a_2\alpha_2^n+\cdots+a_d\alpha_d^n)^{s/m}\notag\\
&=a_1^{s/m}\alpha_1^{ns/m}\left(1+\frac{a_2}{a_1}(\frac{\alpha_2}{\alpha_1})^n+\cdots+\frac{a_d}{a_1}(\frac{\alpha_d}{\alpha_1})^n\right)^{s/m}\notag\\
&=a_1^{s/m}\alpha_1^{ns/m}\sum\limits_{\bar{h}}b_{\bar{h}}a_2^{h_2}\cdots a_d^{h_d}a_1^{-(h_2+\cdots+h_d)}\left(\alpha_2^{h_2}\cdots\alpha_d^{h_d}
\alpha_1^{-(h_2+\cdots+h_d)}\right)^n,
\end{split}
\end{equation}
where $\bar{h}=(h_2,\ldots,h_d)$ runs through $\N^{d-1}$.
Now, we put $y=1/x$. Then we see that $h(x)=h(1/y)$ can be written as an infinite sum
\[
h(1/y)=\uu_{-1}G_n(1/y)^{1/m}+\uu_0+\uu_{1}G_n(1/y)^{-1/m}+\uu_{2}G_n(1/y)^{-2/m}+\cdots
\]
of terms $t(x)=t_{h_2,\ldots,h_d,s}(1/y)$ of the shape
\begin{align}\label{shapeterms}
c\cdot\alpha_1(1/y)^{ns/m-n(h_2+\cdots+h_d)}\cdot\alpha_2(1/y)^{nh_2}\cdots\alpha_d(1/y)^{nh_d},\quad c\in\C.
\end{align}
Observe that as an element of $\C(\!(y)\!)$ the rational function $h(1/y)$ equals $y^{-\deg h}+\cdots$, since by Lemma \ref{lem1} it can be written as $y^{-\deg h}$ times a polynomial in $y$ starting with a non-zero constant term (which is a unit in the ring $\C[\![y]\!]$). (Here we assume, as we may, that $h$ is monic; otherwise the series would start with $y^{-\deg h}$ times the leading coefficient of $h$.) A similar consideration shows that each of the terms on the right hand side, written down explicitly in (\ref{shapeterms}), up to a non-zero constant equals $y^{(-\deg\alpha_1(s/m-(h_2+\cdots+h_d)-h_2\deg\alpha_2-\cdots-h_d\deg\alpha_d)n}$ times a power series that is a unit in the ring $\C[\![y]\!]$. The ``smallest'' such term (i.e. the term with smallest order as a Laurent series) appears precisely with $s=1$ and $h_2=\cdots=h_d=0$ (since by assumption $\deg\alpha_1>\deg\alpha_i$ for $i\ge2$), from which we see that the right hand side starts with $y^{-n\deg\alpha_1/m}$ up to a constant in $\C$. In view of $\deg G_n=n\deg\alpha_1=m\deg h=\deg g\deg h=\deg (g\circ h)$ this perfectly makes sense. We also remark that only terms with $s=0,1$ contribute to $h(1/y)=y^{-\deg h}+\cdots+h_0$, where $h_0=0=h(0)$, since terms with $s<0$, up to a constant, start with $y^{|s|\deg\alpha_1n/m}$. We still have to show that the infinite sum converges as an element of $\C(\!(y)\!)$. This is clear since for $s\rightarrow -\infty$ the power of $y$ goes to $\infty$. We give an alternative proof: We show that for an arbitrary $J\in\N$, there is an upper bound on the number $L$ of quantities $t_1(1/y),\ldots,t_L(1/y)$ in the expansion of $h(1/y)$ that satisfy $\nu_0(t_i(1/y))<nJ$, where $\nu_0$ denotes the order on $\C(\!(y)\!)$ (observe that $n$ is considered at this point to be fixed such that $G_n(x)=g\circ h$). For a term of the shape \eqref{shapeterms}, its order is given by
\begin{align}
\nu_0(t(1/y))&=n\cdot[(\frac{s}{m}-(h_2+\cdots+h_t))\cdot\nu_0(\alpha_1(1/y))+h_2\nu_0(\alpha_2(1/y))+\cdots\notag\\
&\qquad\qquad\qquad\cdots+h_d\nu_0(\alpha_d(1/y))]\notag\\
&=n\cdot[\frac{s}{m}\cdot\nu_0(\alpha_1(1/y))+h_2(\nu_0(\alpha_2(1/y)-\nu_0(\alpha_1(1/y)))+\cdots\notag\\
&\qquad\qquad\qquad\cdots+h_d(\nu_0(\alpha_d(1/y))-\nu_0(\alpha_1(1/y)))].
\end{align}
Note again that by Lemma \ref{lem1} we have $\nu_0(\alpha_i(1/y))=\nu_0(\alpha_i(x))=-\deg\alpha_i$. Since by assumption $\deg\alpha_1>\deg\alpha_i$ for $i\ge2$, it follows that $\deg\alpha_1\ge\deg A_0/d$. In the case $s\le0$ we therefore find
\begin{align*}
\nu_0(t(1/y))&=n[\frac{-s}{m}\deg\alpha_1+h_2(\deg\alpha_1-\deg\alpha_2)+\cdots+h_d(\deg\alpha_1-\deg\alpha_d)]\\
&\ge n[\frac{-s}{m}\cdot\frac{\deg A_0}{d}+h_2+\cdots+h_d].
\end{align*}

We observe that if $\nu_0(t_i(1/y))<nJ$, we must have that $s\in\{0,-1,\ldots,-Jmd/\deg A_0+1\}$, since otherwise we would get
\begin{align*}
\nu_0(t(1/y))\ge n[\frac{-s}{m}\cdot\frac{\deg A_0}{d}+h_2+\cdots+h_d]\ge n[\frac{Jmd}{\deg A_0}\cdot\frac{\deg A_0}{md}]= nJ.
\end{align*}
To estimate the number of possible $(d-1)$-tuples $(h_2,\ldots,h_d)\in\N^{d-1}$ for each $s\in\{0,-1,\ldots,-Jmd/\deg A_0+1\}$, it obviously must hold that \[h_i<J-\frac{|s|\deg A_0}{md}\le J,\] so for each such $s$ there are at most $J^{d-1}$ such $(h_2,\ldots,h_d)$.

If $s=1$ we have that
\begin{align*}
\nu_0(t(1/y))&=n[\frac{-\deg\alpha_1}{m}+h_2(\deg\alpha_1-\deg\alpha_2)+\cdots+h_d(\deg\alpha_1-\deg\alpha_d)]\\
&\ge n[\frac{-\deg\alpha_1}{m}+h_2+\cdots+h_d]\ge n[\frac{-\deg A_0}{m}+h_2+\cdots+h_d],
\end{align*}
hence, if $s=1$, we must have $h_i<J+\frac{\deg A_0}{m}$. Then the number of possible $(d-1)$-tuples $(h_2,\ldots,h_d)$ is not greater than $(J+\frac{\deg A_0}{m})^{d-1}$.

We conclude that we may write
\begin{align}
h(x)=h(1/y)=t_1(1/y)+\cdots+t_L(1/y)+\sum\limits_{\nu_0(t(1/y))\ge Jn}t(1/y),
\end{align}
where $L$ is bounded above by
\begin{align}
L\le \left(J+\frac{\deg A_0}{m}\right)^{d-1}+J^d\frac{md}{\deg A_0}.
\end{align}
This now justifies the above formal expansions. They are well-defined in the ring $\C(\!(y)\!)$.

We now distinguish between two cases, namely that $\{t_1,\ldots,t_L,h(x)\}$ is linearly dependent or linearly independent over $\C$, respectively (we will often simply write $t_i$ instead of $t_i(x)$).\\

\noindent\emph{Case 1. Let us assume that the set $\{t_1,\ldots,t_L,h(x)\}$ is linearly independent over $\C$.}\\

With the intention of applying Proposition \ref{zannier}, let $F=\C(y,\alpha_1(1/y)^{1/m})$ and write $\varphi_1=-t_1(1/y),\ldots,\varphi_L=-t_L(1/y)$ and $\varphi_{L+1}=h(1/y)$. Also, set $\sigma=\sum_{i=1}^{L+1}\varphi_i=\sum_{\nu_0(t(1/y))\ge Jn}t(1/y)$. Observe that $\varphi(T)=T^{m_0}-\alpha_1(1/y)^{m_0/m}$ is the minimal polynomial of $\alpha_1(1/y)^{1/m}$ over $\C(y)$, where we use that the definition of $m_0$ for $\alpha_1$ as a polynomial in $x$ over $\C$ implies that $\alpha_1(1/y)^{m_0/m}$ is not a power of an element in $\C(y)$ for a smaller power. Since this is a Kummer extension we can apply Proposition \ref{kummer} (Theorem III.7.3 of \cite{St}). Therefore, we get that only places in $F$ above $0,\infty$ and the inverses of non-zero roots of $\alpha_1$ (as a polynomial in $\C[x]$) ramify. Thus for the genus $\g_F$ of $F$ we find $2\g_F-2\leq m_0\deg\alpha_1\leq m\deg\alpha_1$. We define $S$ to be the set of zeros and poles of the $t_1,\ldots,t_L$ together with the poles of $h(1/y)$. Observe that $h(1/y)$ has poles at most at places above $0,\infty$. Therefore $S$ may contain at most the places above $0,\infty$ and the inverses of the non-zero roots of $\alpha_1,\ldots,\alpha_d$. This gives at most $m_0(2+\deg\alpha_1+\ldots+\deg\alpha_d)\leq m(2+\deg A_0)$ elements in $S$.
%\begin{align}\label{S}
%S=\{a\in S_0&|\ \nu_a(\alpha_i(y))\neq0\textup{ for some }1\le i\le t\}\cup\{\infty\}\notag\\
%&\cup\{a\in S_0|\ \nu_a(h(y))<0\},
%\end{align}
%where $S_0$ denotes the set of finite places of $F$.

Note that for any place $P$ in $F$ above $0$ in $\C(y)$ we have that
$\nu_P(\sigma)=\nu_P(\sum_{\nu_0(t(1/y))\ge Jn}t(1/y))=e(P|0)\cdot\nu_0(\sum_{\nu_0(t(1/y))\ge Jn}t(1/y))\ge Jn.$
Clearly, $P\in S$.
%By Eisenstein's criterion (see e.g. Theorem III.1.14 in \cite{St}) there is exactly one place $P$ in $F$ above $0$ in $\C(y)$ which totally ramifies and for which $\nu_P(\alpha_1(1/y))^{1/m})=-\deg\alpha_1m_0/m$. Observe that gcd$(\nu_P(\alpha_1(1/y)^{1/m},m_0)=1$ as can be seen by the following argument: Write $\alpha_1(x)=v(x-v_1)^{k_1}\cdots(x-v_t)^{k_t}$. Put $d=$ gcd$(k_1,\ldots,k_t,m)$. Then $m_0=m/d$. It follows that a divisor $d$ of $m_0$ is a divisor of $m$ which cannot be a common divisor of $k_1,\ldots,k_t$. Hence unless $d$ is trivial, we get a contradiction to the definition of $m_0$. Therefore, the conclusion about the order of vanishing of $\sigma$ translates to $\nu_P(\sigma)\geq m_0\nu_0(\sigma)\geq Jnm_0$. Clearly, $P\in S$.
%Since $F$ is just the rational function field, we find $\g=0$ for the genus. To estimate $|S|$, note that $\alpha_i(y)$ has exactly one pole (at $x=0$) and just the same number of zeros as $\alpha_i(x)$ or possibly one less (if $\alpha_i(0)=0$). Since $A_0(x)=\prod_{i=1}^t\alpha_i(x)$ by Vieta, we find that
%\begin{align}
%\#\{a\in S_0|&\ \nu_a(\alpha_i(y))\neq0\textup{ for some }1\le i\le t\}\notag\\
%&\le 1+\#\{a\in S_0|\ \nu_a(\alpha_i(x))\neq0\textup{ for some }1\le i\le t\}\notag\\
%&= 1+\#\{a\in S_0|\ \nu_a(A_0(x))\neq0\}\le1+\deg A_0.
%\end{align}
%The only pole of $h(y)$ is at $x=0$. Since we have already taken this place into account, we get that
%\[|S|\le\deg A_0+2.\]
We will also need to give an upper bound on the degree $\deg h(1/y)=[F:\C(h(1/y))]=\HH(h(1/y))$. Note that $\HH(h(1/y))=(\deg h)\HH(1/y)=(\deg h)[F:\C(1/y)]=(\deg h)[F:\C(x)]=m_0\deg h\leq m\deg h.$ Hence, $\deg h(1/y)= m_0\deg h$.
%from now on we will simply write $\deg h$ instead of $\deg h(x)=\deg h(y)$.
%Furthermore, by the recurrence relation \eqref{recrel} we inductively find the following upper bound on the degree of $G_n(x)$:
%\begin{align}
%\deg G_n\le(n+1-d)K_1+K_2\le n\cdot K,
%\end{align}
%where $K_1=\max_{i=0,\ldots,d-1}(\deg A_i)$, $K_2=\max_{i=0,\ldots,d-1}(\deg G_i)$ and $K=\max(K_1,K_2)$.
%It follows from $G_n(x)=g(h(x))$ and $\deg g\ge2$ that
%\[\deg h(y)=\deg h(x)\le nK/2.\]

By Proposition \ref{zannier} we find that
\begin{align}
\sum\limits_{\nu\in S}(\nu(\sigma)-\min_{i=1,\ldots,L+1} \nu(\varphi_i))&\le\frac{1}{2}L(L+1)(|S|+2\g_F-2)+\deg h(1/y)\notag\\
&\le\frac{1}{2}L(L+1)2m_0(\deg A_0+1)+m_0\deg h.
\end{align}

On the other hand, since $\sigma=\sum_{i=1}^{L+1}\varphi_i$, it follows that $\nu(\sigma)-\min_{i=1,\ldots,L+1} \nu(\varphi_i)\ge0$ for every valuation $\nu\in S$. Moreover, for $\nu_P(\sigma)\ge J n$ and $\min_{i=1,\ldots,L+1} \nu_P(\varphi_i)\le \nu_P(h(1/y))\le m_0\nu_0(h(1/y)))=-m_0\deg h$ (by Lemma \ref{lem1}), we see that
\[J n+m_0\deg h\le\nu_P(\sigma)-\min_{i=1,\ldots,L+1} \nu_P(\varphi_i)\le \sum_{\nu\in S}(\nu(\sigma)-\min_{i=1,\ldots,L+1} \nu(\varphi_i)).\]
We conclude that
\begin{align}
J n+m_0\deg h&\le\sum_{\nu\in S}(\nu(\sigma)-\min_{i=1,\ldots,L+1} \nu(\varphi_i))\\ &\le L(L+1)m_0(\deg A_0+2)+m_0\deg h,\notag
\end{align}
and therefore, for $n$ we get the upper bound
\begin{align}\label{bound1}
n\le m L(L+1)(\deg A_0+2)/J
\end{align}
(recall that for $J\in\N$ one may take any natural number to get an upper bound, and that $L$ is bounded above by a constant depending only on $J$, $m=\deg g$ and the recurrence sequence, but not on $n$).\\

\noindent\emph{Case 2. Let us now consider the second case, namely that the set $\{t_1,\ldots,t_L,h(x)\}$ is linearly dependent over $\C$.}\\

We may assume that the $t_1,\ldots,t_L$ are linearly independent, since otherwise we just group together the terms in question properly. Therefore, in a relation of linear dependency, $h(x)$ must appear and we may write $h(x)$ as a linear combination of $t_1(x),\ldots,t_L(x)$, i.e. there are $w_i\in\C$ such that
\begin{align}\label{ld}
h(x)=\sum\limits_{i=1}^L w_it_i(x).
\end{align}

Recall, that the $t_i$'s are all of the shape \eqref{shapeterms}, where $s$ and the $h_i$'s are elements of a finite set of numbers. After possibly renumbering the terms, we may assume that $w_i\neq0$ exactly for $i=1,\ldots,l\le L$ and we get a power sum representation of $h(x)$ of the shape
\[
h(x)=w_1d_1\delta_1^n+\cdots+w_ld_l\delta_l^n,
\]
where we can control the $\delta_i$'s, since they are elements of a finite set, namely
\[
\delta_i\in\{\alpha_1^{s/m-(h_2+\cdots+h_t)}\alpha_2^{h_2}\cdots\alpha_d^{h_d};\ s\in A,\ h_i\in B\},
%d_i&\in\{a_1^{s/m-(h_2+\cdots+h_t)}a_2^{h_2}\cdots a_d^{h_d}|\ s\in A,\ h_i\in B\},
\]
where $A=\{0,1\}$ and $B=\{0,1,\ldots,J+\deg A_0/m\}$ (for $A$ only terms $t_i$ with $s=0,1$ contribute, as we have already observed above). However, note that we have no control over the coefficients $w_id_i\in\C$.

We pause a moment to investigate this relation. Remember that $m_0$ was defined to be the least integer such that $\alpha_1^{m_0/m}\in\C(x)$. Thus we have $m_0=[F:\C(x)]$. Observe that $m_0$ is a divisor of $m$. We may then write $h$ in the following form \[h(x)=\sum_{j=0}^{m_0-1}\alpha_1^{j/m}\Lambda_j,\] where $\alpha_1^{j/m}\Lambda_j$ is the sum of the terms of the shape \eqref{shapeterms} for which $sn\equiv j$ (mod $m_0$); in particular, $\Lambda_j\in\C(x)$. Since $h\in\C[x]$, we deduce that $h(x)=\Lambda_0$. %Taking the trace from $F$ to $\C(x)$ we see that only those summands contribute for which $\delta_i^n\in\C(x)$.
Note that at least one $t_i$ with $s=1$ has to appear since otherwise $h$ would be constant. Moreover, for $G_n(x)=g(h(x))$ to be true, the term with $s=1$ and $h_2,\ldots,h_d=0$ must appear, hence from the special shape of $t_i$ it follows that we therefore necessarily have $n\equiv 0$ (mod $m_0$). We shall write $n=m_0\ell$ from now on.

Putting $w_id_i=c_i$ and $\delta_i^{m_0}=\gamma_i$ for $i=1,\ldots,l$ gives $h(x)=c_1\gamma_1^\ell+\cdots+c_l\gamma_l^\ell$, which is the claim in Theorem \ref{th1}.
\hfill$\square$

\subsection{Proof of Theorem \ref{th2}}

Define $H_\ell=c_1\gamma_1^\ell+\cdots+c_l\gamma_l^\ell\in\C[c_1,\ldots,c_l](x)$ for any $\ell\in\N$, where $l$ and the $\gamma_i$ are determined in Theorem \ref{th1}. Moreover, we set
\[g(x)=g_0x^m+g_{1}x^{m-1}+\cdots+g_m\in\C[g_0,\ldots,g_m][x].\]
We may write $g\circ h$ as a finite sum of $\ell$-th-power terms:
\[
g(h(x))=g(c_1\gamma_1^\ell+\cdots +c_l\gamma_l^\ell)=e_1\epsilon_1^\ell+\cdots+e_k\epsilon_k^\ell.
\]
Here, we can control the $\epsilon_i$'s and $k$, but not the $e_i$'s since they depend on the $w_i$'s in \eqref{ld}. In fact, we have that
\[
%\delta_i&\in\{g_kd_{i_1}\cdots d_{i_k};\ 0\le k\le m,\ 1\le i_j\le l,\ 1\le j\le k\},\label{delt}\\
\epsilon_i\in\{\gamma_{i_1}\cdots\gamma_{i_r};\ 0\le r\le m,\ 1\le i_j\le l,\ 1\le j\le r\},\label{eps}
\]
and
\[k\le 1+l+l^2+\cdots+l^m=(l^{m+1}-1)/(l-1),
\]
and the $e_i$ are polynomial expressions in $g_0,\ldots,g_m,c_1,\ldots,c_l$ with coefficients in $\Q$, which can be written down explicitly. Moreover, we can assume that $\epsilon_i/\epsilon_j$ is not in $\C$ for any $1\leq i<j\leq k$ because otherwise we can join the two terms with the cost that the $e_i$ are polynomial-exponential equations in $g_0,\ldots,g_m,c_1,\ldots,c_l$ and exponentials in $\ell$ with base in $\C$, which again can be written down explicitly. Below we will have that all the $\epsilon_i$'s are polynomials over $\C$ in which case we can put out the leading coefficient, which is another exponential expression in $\ell$ with base in $\C$; we can put these exponentials also inside the $e_i$'s by modifying the $e_i$'s and $\epsilon_i$'s accordingly.

Now from the proof of Theorem \ref{th1} we know that if $G_n(x)=g(h(x))$, then either $n$ is bounded above by \eqref{bound1}  or $n=m_0\ell$ for some $\ell\in\N$ and $h(x)$ may be written as a power sum as above. In this case we get the equation
\begin{align}\label{case2equ}
a_1b_1^{m_0\ell}\beta_1^{m_0\ell}+\cdots+a_db_d^{m_0\ell}\beta_d^{m_0\ell}=G_{m_0\ell}(x)=e_1\epsilon_1^{\ell}+\cdots+e_k\epsilon_k^{\ell}.
\end{align}
We seek to apply Theorem \ref{bm} to this equation to, again, give an upper bound on possible indices $n\in\N$ in this case, if possible. The last equation can be viewed as a homogeneous $S$-unit-equation over $K=\C(x)$.
%We write the equation in the form
%\begin{align}\label{case2eq2}
%1+\frac{\pi_1}{\pi_t}\left(\frac{\alpha_1}{\alpha_t}\right)^n+\ldots+\frac{\pi_{t-1}}{\pi_t}\left(\frac{\alpha_{t-1}}{\alpha_t}\right)^n+\frac{-\delta_1}{\pi_t}\left(\frac{\epsilon_1}{\alpha_t}\right)^n+\ldots+\frac{-\delta_k}{\pi_t}\left(\frac{\epsilon_k}{\alpha_t}\right)^n=0.
%\end{align}
Application of Theorem \ref{bm} requires the equation in question not to have a proper vanishing subsum. Therefore we look at a minimal vanishing subsum of \eqref{case2equ} and again we distinguish between two possible cases. If there is such a subsum consisting of at least three terms, we are able to apply the theorem and consequently get an upper bound on $n$. On the other hand, this is not the case if and only if each of the terms of the right hand side of equation \eqref{case2equ} is identical to  exactly one term on the other side (note that in this case each minimal vanishing subsum consists of exactly two terms, since it cannot consist of one term by minimality, and moreover, because of the non-degeneracy it is not possible that two terms on the same side of equation \eqref{case2equ} coincide). Also, if we are in this exceptional situation, it follows that $k=d$ and that there is a unique permutation $\rho\in \mathscr{S}_d$ such that $a_i\alpha_i^{m_0\ell}=e_{\rho(i)}\epsilon_{\rho(i)}^\ell$ (if there was another permutation of this kind, we would again end up in the situation that two terms of the same side coincide).

Now, assume that there is a proper vanishing subsum of \eqref{case2equ} consisting of at least three terms. Then we may write the equation in the form $1+v_1+\cdots+v_n=0$ (where $2\le n\le k+d-1$). Since for zeros and poles of the $v_i$'s we can only have zeros and poles of the $\alpha_i$'s, we can take $K=\C(x)$ and $S_K$ the set of places of $K$ containing $\infty$ and the zeros of $\alpha_1,\ldots,\alpha_d$. Similarly as above it follows that $|S_K|\le \deg A_0+1$ and clearly $\g_K=0$,
%\begin{align*}
%|S|\le \deg A_0+2\quad\text{ and }\quad\g=0,
%\end{align*}
hence by Brownawell and Masser's theorem (Theorem \ref{bm}) it follows that
\begin{align}\label{upper}
\max\limits_{i=1,\ldots,n}\HH(v_i)&\le\frac{1}{2}(k+d-2)(k+d-3)(|S_K|+2\g_K-2)\notag\\
&\le \frac{1}{2}(\frac{L^{m+1}-1}{L-1}+d-2)^2\deg A_0.
\end{align}
Assume first that the vanishing subsum contains $\beta_i$ and $\beta_j$ for $1\leq i<j\leq d$. The case when it contains two of the $\epsilon_i$'s will be done afterwards. Now, we may obtain the equation of shape $1+v_1+\cdots+v_n=0$ by dividing the original one by $a_jb_j^{m_0\ell}\beta_j^{m_0\ell}$. Then we find the following lower bound (note that we have $\HH(\alpha_i/\alpha_j)\ge1$, by the non-degeneracy of the sequence)
\begin{align}\label{lower}
\max\limits_{i=1,\ldots,n}\HH(v_i)&\ge\HH((\alpha_i/\alpha_j)^{m_0\ell})= m_0\ell\cdot \HH(\alpha_j/\alpha_t)\ge m_0\ell=n\ge \ell.
\end{align}
Now assume that the vanishing subsum contains $\epsilon_i$ and $\epsilon_j$ for $1\leq i<j\leq k$. Then as above (note now that we have $\HH(\epsilon_i/\epsilon_j)\geq 1$ because we have taken the $\epsilon_i$'s to be non-degenerate)
\[\max\limits_{i=1,\ldots,n}\HH(v_i)\ge \HH((\epsilon_i/\epsilon_j)^\ell)=\ell\cdot\HH(\epsilon_i/\epsilon_j)\ge \ell.
\]
Clearly, with $\ell$ we also have an upper bound for $n=m_0\ell$. Hence, by \eqref{upper} and \eqref{lower} we see, that again in this case there can only be finitely many $n\in\N$ such that $G_n(x)=g(h(x))$ for some polynomial $h\in\C[x]$, with an upper bound $C'$ given by \eqref{upper}.

As we saw, the only remaining situation we have to consider is the situation that in equation \eqref{case2equ} each term on the left hand side coincides with exactly one term on the right hand side. In particular, it holds that $d=k$. So let $\rho\in \mathscr{S}_d$ be such that $a_i\alpha_i^{m_0\ell}=a_ib_i^{m_0\ell}\beta_i^{m_0\ell}=e_{\rho(i)}\epsilon_{\rho(i)}^\ell$ for all $i=1,\ldots,d$, or equivalently
\begin{equation}\label{exceptional}
\left(\frac{\beta_i^{m_0}}{\epsilon_{\rho(i)}}\right)^\ell=\frac{e_{\rho(i)}}{a_ib_i^{m_0\ell}}
\end{equation}
for all $i=1,\ldots,d$. Since $a_i,b_i,e_{\rho(i)}\in\C$ are constant, it follows that $\beta_i^{m_0}$ and $\epsilon_{\rho(i)}$ are polynomials that coincide up to some constant factor. We have already mentioned above that we may assume that both polynomials have leading coefficient equal to one. Therefore it follows that they have to be equal so that the quotient is equal to one for all $\ell$. It follows that we have equalities $a_ib_i^{m_0\ell}=e_{\rho(i)}$ for $i=1,\ldots,d$. These are polynomial-exponential equations in the unknowns $g_0,\ldots,g_m,c_1,\ldots,c_l$ and $\ell$ which define a subvariety as claimed in the theorem.

We have shown that there is a subvariety $\mathcal{V}$ of $\mathbb{A}^{l+m+1}\times\mathbb{G}_\m^{t}$, where $l$ and $t$ are explicitly bounded by the originally given data and which is given by polynomial-exponential equations in the polynomial unknowns $g_0,\ldots,g_m,c_1,\ldots,c_l$ and the exponential unknown $\ell$ which can be written down explicitly. If we then define $G(x)=g_0x^m+\cdots+g_m\in\C[\mathcal{V}][x]$ and $H_\ell=c_1d_1\delta_1^\ell+\cdots+c_ld_l\delta_l^\ell\in\C[\mathcal{V}](x)$, then by construction we have for all $\ell\in\N$ that $G_{m_0\ell}(x)=G\circ H_\ell$ as an equation in $\C[\mathcal{V}](x)$ and the following holds: If $G_n(x)=g\circ h$ with $g,h\in\C[x],\deg g=m,\deg h>1$, then either $n$ is not greater than the maximum of the bounds given in \eqref{bound1} and $C'$ or $n=m_0\ell$ and there is a point $P\in\mathcal{V}(\C)$ such that $g(x)=G(P,x)$ and $h(x)=H_\ell(P,x)$. Conversely, if $P=(g_0,\ldots,g_m,c_1,\ldots,c_l,\ell)\in\mathcal{V}(\C)$ is a given $\C$-rational point on $\mathcal{V}$ and if we define $g(x)=G(P,x)$ and $h(x)=H_\ell(P,x)$ then we have $G_{m_0\ell}(x)=g\circ h$ with $\deg g=m$. By integrality it follows that $h\in\C[x]$. Clearly, $\deg h>1$. This establishes the theorem.
\hfill $\square$

\section{Proof of the remarks and special cases}
\label{sec4}

a) We just have to apply \cite[Theorem 1]{FKK}. It follows that there is a constant $C$ with $\deg g\leq C$. We apply our Theorem \ref{th1} for each $m$ with $2\leq m\leq C$. From this the conclusion follows.

b) This follows directly by inspecting the proof. The equations for the varieties arise from (\ref{exceptional}). By assumption $b_i=1$, $\beta_i,\epsilon_j$ are monic and $e_j$ is  a polynomial in $c_1,\ldots,c_l,g_0,\ldots,g_m$ with rational coefficients. We therefore get for $\mathcal{V}$ a subvariety of $\mathbb{A}^{l+m+1}$ and there is no $\mathbb{G}_{\textnormal{m}}$-part.

c) Assume that $G_n(x)=g\circ h$ with $\deg g=m$ and $n>C$. Then by Theorem \ref{th1} it follows that $n=m_0\ell$ and that $h(x)=c_1\gamma_1^\ell+\cdots+c_l\gamma_l^\ell$, where the $\gamma_i$ are (up to a constant) of the shape
\[\left(\alpha_1^{s_i/m-(h_{i2}+\cdots+h_{id})}\alpha_2^{h_{i2}}\cdots\alpha_d^{h_{id}}\right)^{m_0},\] where $s_i\in A=\{0,1\},h_{ij}\in B=\{0,1,\ldots,J+\deg A_0/m\}$. We put $k=\deg\alpha_1$ and assume that $1\le k< m$.
In the case $(h_{i2},\ldots,h_{id})=(0,\ldots,0)$, if $s=0$ then  $\gamma_i$ is constant, and if $s=1$ then $\gamma_i=\xi\alpha_1^{m_0/m}\in\C[x]$, where $\xi\in\C$.
Assume that there is an $i$ with $(h_{i2},\ldots,h_{id})\neq(0,\ldots,0)$ and let ${\bf{h}}=\max_i\sum_{j=2}^d h_{ij}$.
Then $\gamma_i$ (up to a constant) is of the shape
\begin{align*}
\left(\frac{\alpha_1^{s_i/m-(h_{i2}+\cdots+h_{id})+{\bf{h}}}\alpha_2^{h_{i2}}\cdots\alpha_d^{h_{id}}}{\alpha_1^{{\bf{h}}}}\right)^{m_0}
=\left(\frac{p(x)}{\alpha_1^{{\bf{h}}}}\right)^{m_0},
\end{align*}
where $p(x)\in\C[x]$ is a polynomial with
\begin{align*}
\deg p&=ks_i/m-k(h_{i2}+\cdots+h_{id})+{\bf{h}}k+h_{i2}\deg\alpha_2+\cdots+ h_{id}\deg\alpha_d\\
&<1+{\bf{h}}k+\sum_{j=2}^d h_{ij}(\deg\alpha_j-k)\le {\bf{h}}k.
\end{align*}
%Therefore, whenever $(h_{i2},\ldots,h_{id})\neq(0,\ldots,0)$, then $\gamma_i$ is of the shape $(p(x)/x^{{\bf{h}}k})^{m_0}$, where $\deg p<{\bf{h}}k$.
Hence, subtracting possibly the terms with $(h_2,\ldots,h_d)=(0,\ldots,0)$ from $h(x)$, we get a polynomial $h'(x)$ which can be written as the sum of $\ell$-th powers of terms of the described shape, i.e. $h'(x)=p'(x)/\alpha_1^{{\bf{h}}m_0\ell}$, where $\deg p'<{\bf{h}}k m_0\ell$. This gives a contradiction to $h'\in\C[x]$, except for the case $p'=0$.
Thus we may assume that only the summands with $h_2=\cdots=h_t=0$ (and $s=0,1$) occur, that is $h(x)=c_1+c_2\alpha_1^{(m_0\ell)/m}=c_1+c_2\alpha_1^{n/m}=c_1+c_2\alpha_1^{\deg G_n/(mk)}$ (since $k=\deg\alpha_1>\deg\alpha_i$ for $i=2,\ldots,d$ we have that $\deg G_n=n\deg\alpha_1=kn$). A similar argument (considering the valuation at 0) shows that $h(x)$ must also be of this shape in the case $\deg\alpha_1=k=m$. Note that for the finitely many remaining cases when $m<k$ one can apply Theorems \ref{th1} and \ref{th2}, respectively.

%in the function field $\C(x)$. By assumption $\nu$ is not a zero (and certainly not a pole) of $\alpha_2,\ldots,\alpha_t$. Assume that there is an $i$ with $(h_{i2},\ldots,h_{id})\neq(0,\ldots,0)$. Then $\nu(\gamma_i)=s/m-(h_{i2}+\cdots+h_{id})<0$. It follows by the strict triangle inequality that $0\leq \nu(h)=\nu(\sum c_j\gamma_j^\ell)$, where the sum is taken over all summands $c_j\gamma_j^\ell$ with $\nu(\gamma_j)$ minimal. Since each summand of $\sum c_j\gamma_j^\ell$ has a pole at $\nu$ and since by the special shape of the summands the sum has also a pole at $\nu$ (here we use that $\deg\alpha_i<\deg\alpha_1$ for $i=2,\ldots,t$) we get a contradiction. Thus only the summands with $h_2=\cdots=h_t=0$ (and $s=0,1$) occur proving the statement.

d) We argue as in c). We therefore have that $\gamma_i=\beta^{(s_i/m-(h_{i2}+\cdots+h_{id}))m_1+h_2m_2+\cdots+h_dm_d}=\beta^{k_i}$ for some $k_i\in\Z$. By integrality it follows that $k_i\in\N$. This shows the claim.

e) We assume that $a_1,\ldots,a_d\in K$ and $\alpha_1,\ldots,\alpha_d\in K[x]$ so that $G_n(x)\in K[x]$ for all $n\geq 0$, where $K$ is a number field. An application of Theorem \ref{th1} implies that either $n\leq C$ or $h(x)=c_1\gamma_1^\ell+\cdots+c_l\gamma_l^\ell$ with $n=m_0\ell$. We have $\gamma_1,\ldots,\gamma_l\in K(x)$ and $c_1,\ldots,c_l\in\C$. Taking a transcendence basis of $K(c_1,\ldots,c_l)\supseteq K$ and comparing coefficients in this basis shows that we may assume that $c_1,\ldots,c_l$ are algebraic over $K$. Taking now a (field) basis of $K(c_1,\ldots,c_l)\supseteq K$ and comparing coefficients in this basis shows that we may assume $c_1,\ldots,c_l\in K$. Since $G_n(x)\in K[x]$ and $h(x)\in K[x]$ determine $g$, we see that $g(x)\in K[x]$ as well (for this the algorithm in \cite[Secion 3]{R} can be used). This proves the claim.

f) This follows immediately by fixing $g(x)=x^m$ and then going through the proofs.

g) Assume that $h(x)=c_1\gamma_1^\ell+\cdots+c_l\gamma_l^\ell$ with $n=m_0\ell$ appears as ``inner'' decomposition factor for fixed $c_1,\ldots,c_l$ (independent of $\ell$) for any $\ell\in\N$. Then, as already mentioned above, we can identify the varying powers by indeterminates and then use the algorithm provided in \cite{B} to explicitly calculate this family of decompositions.

%\newpage
%\section*{Acknowledgement} The author was supported by Austrian Science Fund (FWF) Grant No. P24574.


\begin{thebibliography}{XX}
\bibitem{AZ} {\sc R.M.\  Avanzi and U.\ Zannier}, The equation $f(X)=f(Y)$  in rational functions $X=X(t)$ , $Y=Y(t)$. Compositio Math.  139  (2003),  no. 3, 263--295.
\bibitem{BT} {\sc Y.\ Bilu and R.\ Tichy}, The Diophantine equation $f(x)=g(y)$. Acta Arith.  95  (2000),  no. 3, 261--288.
\bibitem{Bea}{\sc  Y.\ Bilu, C.\ Fuchs, F.\ Luca, and \'{A}.\ Pint\'{e}r}, Combinatorial Diophantine equations and a refinement of a theorem on separated variables equations. Publ. Math. Debrecen  82  (2013),  no. 1, 219--254.
\bibitem{B} {\sc A.\ Bodin}, Decomposition of polynomials and approximate roots. Proc. Amer. Math. Soc.  138  (2010),  no. 6, 1989--1994.
\bibitem{BM}{\sc W.D.\ Brownawell and D.W.\ Masser}, Vanishing sums in function fields. Math. Proc. Cambridge Philos. Soc. 100 (1986), no. 3, 427--434.
\bibitem{CZ} {\sc P.\ Corvaja and U.\ Zannier}, Finiteness of integral values for the ratio of two linear recurrences. Invent. Math.  149  (2002),  no. 2, 431--451.
%\bibitem{CZ} Corvaja, Pietro; Zannier, Umberto Some new applications of the subspace theorem. Compositio Math.  131  (2002),  no. 3, 319–340.
%\bibitem{CZ2} Corvaja, Pietro; Zannier, Umberto S -unit points on analytic hypersurfaces. Ann. Sci. École Norm. Sup. (4)  38  (2005),  no. 1, 76–92.
\bibitem{E}{\sc M. Eichler}, Einf\"uhrung in die Theorie der algebraischen Zahlen und Funktionen. Lehrb\"ucher und Monographien aus dem Gebiete der exakten Wissenschaften, Mathematische Reihe, Nr. 27, Birkh\"auser Verlag, Basel 1963.
%\bibitem{E2} {\sc J.-H.\ Evertse}, Unit equations in Diophantine number theory. Cambridge University Press, Cambridge, 2015.
\bibitem{F} {\sc C. Fuchs}, On the Diophantine equation $G_n (x)=G_m (P(x))$  for third order linear recurring sequences. Port. Math. (N.S.)  61  (2004),  no. 1, 1--24.
\bibitem{FKK} {\sc C. Fuchs, C. Karolus, and D. Kreso} Decomposable polynomials in second order linear recurrence sequences. manuscripta math. (https://doi.org/10.1007/s00229-018-1070-8).
\bibitem{FMZ} {\sc C.\ Fuchs, V.\ Mantova, and U.\ Zannier} On fewnomials, integral points, and a toric version of Bertini's theorem. J. Amer. Math. Soc.  31  (2018),  no. 1, 107--134.
%\bibitem{FS} Fuchs, Clemens; Scremin, Amedeo Polynomial-exponential equations involving several linear recurrences. Publ. Math. Debrecen  65  (2004),  no. 1-2, 149–172.
\bibitem{FP1} {\sc C.\ Fuchs and A.\ Peth\H{o} }, Effective bounds for the zeros of linear recurrences in function fields. J. Th\'eor. nombres Bordeaux 17 (2005), 749--766.
\bibitem{FP2} {\sc C.\ Fuchs and  A.\ Peth\H{o}}, Composite rational functions having a bounded number of zeros and poles. Proc. Amer. Math. Soc.  139  (2011),  no. 1, 31--38
\bibitem{FPT1}{\sc C.\ Fuchs,  A.\ Peth\H{o}, and R.F. Tichy}, On the Diophantine equation $G_n(x)=G_m(P(x))$. Monatsh. Math.  137  (2002),  no. 3, 173--196.
\bibitem{FPT2} {\sc C. \ Fuchs,  A.\ Peth\H{o}, and R.F.\ Tichy}, On the Diophantine equation $G_n(x)=G_m(P(x))$: higher-order recurrences. Trans. Amer. Math. Soc.  355  (2003),  no. 11, 4657--4681.
\bibitem{FPT3} {\sc C.\ Fuchs,  A.\ Peth\H{o}, and R.F.\ Tichy}, On the Diophantine equation $G_n(x)=G_m(y)$  with $Q(x,y)=0$.  Diophantine approximation,  199–-209, Dev. Math., 16, Springer Wien-New York, Vienna, 2008.
\bibitem{FZ} {\sc C.\ Fuchs and U.\ Zannier}, Composite rational functions expressible with few terms. J. Eur. Math. Soc. (JEMS) 14 (2012), 175--208.
%\bibitem{M}{\sc R.C.\ Mason}, Diophantine equations over function fields. Cambridge University Press, Cambridge et. al., 1984.
%\bibitem{PT} {\sc A.\ Peth\H{o} and S.\ Tengely}, On composite rational functions. {\it Number theory, ana-lysis, and combinatorics}, 241-259, {\it De Gruyter Proc. Math.}, De Gruyter, Berlin, 2014.
\bibitem{R}{\sc  J.\ Rickards}, When is a polynomial a composition of other polynomials? Amer. Math. Monthly  118  (2011),  no. 4, 358--363.
\bibitem{Setal}{\sc H.\ Salzmann, T.\ Grundh\"ofer, H.\ H\"ahl, and R.\ L\"owen}, The classical fields. Structural features of the real and rational numbers. Encyclopedia of Mathematics and its Applications, 112. Cambridge University Press, Cambridge, 2007.
\bibitem{S} {\sc A.\ Schinzel}, Polynomials with special regard to reducibility. With an appendix by Umberto Zannier. Encyclopedia of Mathematics and its Applications, 77. Cambridge University Press, Cambridge, 2000.
\bibitem{SS} {\sc H.P.\ Schlickewei and W.M.\ Schmidt}, The intersection of recurrence sequences. Acta Arith.  72  (1995),  no. 1, 1--44.
\bibitem{St}{\sc H.\ Stichtenoth}, Function Fields and Codes, Universitext, Springer-Verlag, Berlin, 1993.
%\bibitem{za1} {\sc U.\ Zannier}, On the number of terms of a composite polynomial. {\it Acta Arith.} {\bf 127 (2)} (2007), 157-167.
\bibitem{W}{\sc R.J.\ Walker}, Algebraic curves.  Dover Publications, New York, 1962.
\bibitem{Z3} {\sc U. Zannier}, On the number of terms of a composite polynomial. Acta Arith.  127  (2007),  no. 2, 157--167.
\bibitem{Z4} {\sc U. Zannier}, Addendum to the paper: ``On the number of terms of a composite polynomial''. Acta Arith.  140  (2009),  no. 1, 93--99.
\bibitem{Z1} {\sc U.\ Zannier}, On composite lacunary polynomials and the proof of a conjecture of Schinzel. Invent. Math. 174 (2008), no. 1 127--138.
\bibitem{Z2}{\sc U.\ Zannier}, On the integer solutions of exponential equations in function fields. Ann. Inst. Fourier (Grenoble)  54 (2004),  no. 4, 849--874.
\end{thebibliography}
\end{document}